\documentclass[11pt, reqno]{amsart}

\usepackage{amssymb}

\evensidemargin\oddsidemargin

\begin{document}
\baselineskip=18pt
\setcounter{page}{1}

\newcommand{\eqnsection}{
\renewcommand{\theequation}{\thesection.\arabic{equation}}
    \makeatletter
    \csname  @addtoreset\endcsname{equation}{section}
    \makeatother}
\eqnsection
   

\def\a{\alpha}
\def\B{{\bf B}}
\def\cM{{\mathcal{M}}} 
\def\cD{{\mathcal{D}}} 
\def\cG{{\mathcal{G}}} 
\def\cH{{\mathcal{H}}} 
\def\cL{{\mathcal{L}}} 
\def\cI{{\mathcal{I}}}
\def\CC{{\mathbb{C}}} 
\def\Dap{{\rm D}_{0+}^\a} 
\def\Dm{{\rm D}_{-}^\a} 
\def\Dp{{\rm D}_{+}^\a} 
\def\Ea{E_\a}
\def\esp{{\mathbb{E}}} 
\def\F{{\bf F}}
\def\Farl{{\F}_{\a,\lbd,\rho}}
\def\bF{{\bar F}}
\def\bG{{\bar G}}
\def\G{{\bf G}}
\def\Ga{{\Gamma}} 
\def\Gal{{\G_{\a,\lbd}}} 
\def\GG{{\bf \Gamma}}
\def\ii{{\rm i}} 
\def\Iap{{\rm I}_{0+}^\a} 
\def\Im{{\rm I}_{-}^\a} 
\def\Ip{{\rm I}_{+}^\a} 
\def\L{{\bf L}}
\def\lbd{\lambda}
\def\lacc{\left\{}
\def\lcr{\left[}
\def\lpa{\left(}
\def\lva{\left|}
\def\M{{\bf M}}
\def\NN{{\mathbb{N}}} 
\def\pb{{\mathbb{P}}}
\def\ZZ{{\mathbb{Z}}} 
\def\R{{\bf R}}
\def\rl{{\mathbb{R}}}
\def\racc{\right\}}
\def\rpa{\right)}
\def\rcr{\right]}
\def\rva{\right|}
\def\sga{\sigma^{(\a)}}
\def\sgb{\sigma^{(\beta)}}
\def\T{{\bf T}}
\def\Un{{\bf 1}}
\def\X{{\bf X}}
\def\Y{{\bf Y}}
\def\E{{\bf E}}
\def\Z{{\bf Z}}
\def\Warl{{\W}_{\a,\lbd,\rho}}
\def\Za{{\Z_\a}}

\def\elaw{\stackrel{d}{=}}
\def\claw{\stackrel{d}{\longrightarrow}}
\def\elaw{\stackrel{d}{=}}
\def\qed{\hfill$\square$}


\newtheorem{prop}{Proposition}[section]
\newtheorem{coro}[prop]{Corollary}
\newtheorem{lem}{Lemma\!\!}
\newtheorem{ex}[prop]{Example}
\newtheorem{exs}[prop]{Examples}
\newtheorem{rem}{Remarks\!\!\!}
\newtheorem{theo}{Theorem\!\!}
\newtheorem{conj}{Conjecture\!\!}

\renewcommand{\thetheo}{}
\renewcommand{\theconj}{}
\renewcommand{\thelem}{}
\renewcommand{\therem}{}

\newcommand{\noi}{\noindent}
\newcommand{\dis}{\displaystyle }

\title{On the Green function of the killed fractional Laplacian on the periodic domain}

\author[T.~Simon]{Thomas Simon}

\address{Laboratoire Paul Painlev\'e, UMR 8524, Universit\'e de Lille,  Cit\'e Scientifique, F-59655 Villeneuve d'Ascq Cedex, France. {\em Email}: {\tt thomas.simon@univ-lille.fr}}

 


\begin{abstract} 

We give a very simple proof of the positivity and unimodality of the Green function for the killed fractional Laplacian on the periodic domain. The argument relies on the Jacobi triple product and a probabilistic representation of the Green function. We also show by a contour integration that the Green function is completely monotone on the positive part of the periodic domain.

\end{abstract}

\maketitle

\section{Introduction}

\label{Fürth}

In the recent paper \cite{LP}, the trigonometric series
$$G_{\a, c}(x) \; =\; \frac{1}{2\pi} \sum_{n\in \ZZ} \frac{\cos (nx)}{c + \vert n\vert^\a} \; =\; \frac{1}{2\pi} \sum_{n\in \ZZ} \frac{e^{\ii nx}}{c + \vert n\vert^\a}$$
was considered for $\a, c > 0$ as the solution to the periodic boundary value problem  
$$\lcr c + (-\Delta)^{\a/2}\rcr G_{\a, c} (x)\; =\; \delta (x), \qquad x\in [-\pi, \pi],$$
where $(-\Delta)^{\a/2}$ is the fractional Laplacian on the normalized periodic domain $[-\pi,\pi]$ and $\delta$ stands for the Dirac delta distribution. The function $G_{\a, c} (x)$ is called the Green function associated to the periodic operator $ c + (-\Delta)^{\a/2},$ which is for $\a\in (0,2]$ the infinitesimal generator of a symmetric $\a-$stable L\'evy process on the circle killed at an independent exponential time with parameter $c$. More precisely, if $\{Z^{(\a)}_t, \, t\ge 0\}$ denotes a symmetric $\a-$stable L\'evy process on $\rl$ with Fourier transform $\esp[e^{\ii x Z^{(\a)}_t}] = e^{-t\vert x\vert^\a}, \, \tau_c$ an independent random time with density $c e^{-ct}$ on $\rl^+,$ and $\partial$ some cemetery point, the above operator is the infinitesimal generator of the Markov process
$$X^{(\a, c)}_t \; =\; \lacc \begin{array}{ll}
Z^{(\a)}_t \, - \, 2\pi \lcr {\displaystyle \frac{\pi + Z^{(\a)}_t}{2\pi}}\rcr & \mbox{if $t < \tau_c,$}\\

\partial & \mbox{if $t \ge \tau_c.$}\end{array}\right.$$
Throughout, we refer to \cite{ST} for an account on stable processes. Notice that the relationship between $G_{\a, c}(x)$ and the travelling periodic waves of the fractional Korteweg-de Vries equation with speed $c$ is discussed in the introduction of \cite{LP}. See also the references therein for many other related topics. The main result of \cite{LP} is that $G_{\a, c}(x)$ is positive and decreasing on $(0,\pi)$ for every $\a\in (0,2]$ and $c > 0$ - see Theorem 1.1 therein. In this note, we shall obtain the following extended property.  

\begin{theo}
\label{Posuni}
For every $\a\in (0,2]$ and $c > 0,$ the function $G_{\a,c}(x)$ is completely monotone on $(0,\pi).$
\end{theo}   

The proof of Theorem 1.1 in \cite{LP} relies on an integral representation of the Green function with the generalized Mittag-Leffler function $E_{\a,\a} (-cx^\a)$ and, in the case $\a\in (1,2)$, on an analysis involving the Lax-Milgram theorem and the fractional P\'olya-Szeg\"o inequality. In this note we first present an alternative probabilistic proof of Theorem 1.1 in \cite{LP},  based on the stable subordinator and the Jacobi triple product, which is easy and very short. We then show the complete monotonicity, which is first obtained on the half-line for the non-periodic Green function by a contour integration and then transferred to the function $G_{\a,c} (x)$ on $(0,\pi)$ by a summation argument. We also observe in Remark (b) below that $G_{\a, c} (x)$ cannot be extended to a completely monotone function on the half-line, because $G_{\a,c}'(\pi-) = 0.$ In particular, there is no Bernstein representation for $G_{\a,c} (x).$ 
 
\section{Proof of the Theorem}
 
We shall begin with a short and easy proof of Theorem 1.1 in \cite{LP}. Consider the $\beta-$stable subordinator $\{\sgb_t,\, t\ge 0\}$ with $\beta = \a/2\,\in (0,1],$ that is the increasing $\beta-$stable L\'evy process on $\rl^+$ having Laplace transform $\esp[e^{-\lbd \sgb_t}]\; =\; e^{-t\lbd^\beta}.$ Setting $\tau_c$ for an independent exponential random variable with density $c e^{-ct}$ and $\X_{\a, c} = \sigma^{(\beta)}_{\tau_c},$ we have
$$\esp [e^{-n^2 \X_{\a, c}}] \; =\; \frac{c}{c + \vert n\vert^\a}$$
for every $n\in \ZZ.$ Hence, for every $x\in (0,\pi),$ one has 
$$G_{\a, c} (x) \; = \;\frac{1}{2\pi c} \sum_{n\in \ZZ} e^{\ii nx}\, \esp [e^{-n^2 \X_{\a, c}}]\; = \; \frac{1}{2\pi c}\; \esp\lcr \sum_{n\in \ZZ} e^{\ii nx}\, e^{-n^2 \X_{\a, c}}\rcr$$
where the switching between with the expectation and the sum will be justified soon afterwards. The Jacobi triple product - see e.g. Section 10.4 in \cite{AAR} - implies 
\begin{equation}
\label{JTP}
G_{\a, c} (x) \; =\; \frac{1}{2\pi c}\;\esp\lcr \prod_{n\ge 0} \lpa 1 + 2 \cos (x)\, e^{-(2n+1) \X_{\a, c}} + e^{-(4n+2) \X_{\a, c}}\rpa\lpa 1- e^{-(2n+1) \X_{\a, c}}\rpa\rcr
\end{equation}
and shows that $G_{\a, c} (x)$ is positive and decreasing on $(0,\pi),$ since it is the case for the product. It remains to justify the switching, which is plain for $\a > 1$ by Fubini's theorem. The argument for $\a\le 1$ is standard. First, a classical Abel summation argument with the Dirichlet kernel - see e.g. Chapter 1.2 in \cite{Z} - implies, for every $q > p > 0,$
$$\lva\sum_{p \le \vert n\vert \le q} e^{\ii nx}\, e^{-n^2 \X_{\a, c}} \rva\; \le \;\frac{2\, e^{-p^2 \X_{\a, c}}}{\sin(x/2)}\cdot$$ 
Letting $q\to\infty,$ we get
$$\esp\lcr \lva\sum_{\vert n\vert \ge p} e^{\ii nx}\, e^{-n^2 \X_{\a, c}} \rva\rcr \; \le \; \frac{2}{(\sin(x/2))\, (c+ p^\a)}$$
and the same argument implies
$$\lva\sum_{\vert n\vert \ge p} \frac{e^{\ii nx}}{c + \vert n\vert^\a} \rva \; \le \; \frac{2}{(\sin(x/2))\, (c+ p^\a)}\cdot$$
The switching argument follows then from the linearity of the expectation and the triangle inequality, letting $p\to\infty.$ We omit details. \

We now show the complete monotonicity property, which cannot be deduced directly from (\ref{JTP}) since the derivative of the product vanishes at zero. Instead, we will use the summation formula 
\begin{equation}
\label{Sum}
G_{\a, c} (x)\; =\; \sum_{n\in\ZZ} H_{\a, c} (x- 2n\pi), \qquad x\in (0,\pi)
\end{equation}
which is given in (1.6) of \cite{LP}, where $H_{\a, c} (x)$ is the Green function of the operator $c + (-\Delta)^{\a/2}$ on the line. Recall that by Fourier inversion - see e.g. (A.3) in \cite{FL}, one has for every $x > 0$
$$H_{\a, c}(x) \; =\; \frac{1}{2\pi} \int_{\rl} \frac{\cos (tx)}{c + \vert t\vert^\a}\, dt \; =\; \frac{1}{2\pi}\lpa \int_0^{\infty} \frac{e^{\ii tx}}{c + t^\a}\, dt\, +\, \int_0^{\infty} \frac{e^{-\ii tx}}{c + t^\a}\, dt\rpa.$$
Supposing first $\a\in (0,2),$ a contour integration on the first resp. fourth quadrant for the first resp. second integral on the right-hand side implies, after the necessary simplifications,  
\begin{equation}
\label{HC}
H_{\a, c}(x) \; =\; \frac{\sin(\pi\a/2)}{\pi} \int_0^{\infty} e^{-tx}\,\lpa \frac{t^\a}{c^2 + 2c \cos(\pi\a/2)\, t^\a + t^{2\a}}\rpa\, dt
\end{equation}
for every $x > 0.$ This shows that $H_{\a, c} (x)$ is completely monotone on $(0,\infty).$ By parity, we next transform (\ref{Sum}) into 
\begin{equation}
\label{SumB}
G_{\a, c} (x) \; = \; \sum_{n\ge 0} \lpa H_{\a, c} (x+ 2n\pi)\, +\, H_{\a, c} (2(n+1)\pi-x)\rpa
\end{equation}
for every $x\in (0,\pi).$ Setting $F_n (x) = H_{\a, c} (x+ 2n\pi) + H_{\a, c} (2(n+1)\pi-x)$, we compute for every $x \in (0,\pi)$ and $p,n \in\NN$ 
$$F_n^{(2p)} (x) \; = \; \frac{\sin(\pi\a/2)}{\pi} \int_0^{\infty}  \frac{t^{\a + 2p} \lpa e^{-t(x + 2n\pi)} + e^{-t(2(n+1)\pi -x)}\rpa}{c^2 + 2c \cos(\pi\a/2) t^\a + t^{2\a}}\, dt\; > \; 0$$
and
$$F_n^{(2p+1)} (x) \; = \; \frac{\sin(\pi\a/2)}{\pi} \int_0^{\infty}  \frac{t^{\a + 2p+1} \lpa  e^{-t(2(n+1)\pi -x)} - e^{-t(x + 2n\pi)}\rpa}{c^2 + 2c \cos(\pi\a/2) t^\a + t^{2\a}}\, dt\; < \; 0.$$
By (\ref{SumB}) and Fubini's theorem, this shows that 
$$(-1)^p G_{\a, c}^{(p)} (x)\; >\; 0$$ 
for every $x \in (0,\pi)$ and $p\in\NN,$ as required. Finally, for $\a = 2$ a residue computation on the upper half-plane leads to the standard formula $H_{2,c} (x) \, =\, (1/2\sqrt{c}) e^{-\sqrt{c} x},$ and to the closed formula
$$G_{2,c} (x) \; =\; \frac{e^{-\sqrt{c} x} + e^{-\sqrt{c} (2\pi -x)}}{2\sqrt{c} (1- e^{-2\sqrt{c} \pi})}\; =\; \frac{\cosh(\sqrt{c} (\pi -x))}{2\sqrt{c}\, \sinh(\sqrt{c} \pi)}$$ 
on $(0,\pi),$ which was already derived in (A.5) of \cite{LP} by different means. Clearly, this function is completely monotone on $(0,\pi).$

\qed

\begin{rem} {\em (a) For $\a \in (0,1],$ if we consider the $\a-$stable subordinator $\{\sga_t,\, t\ge 0\}$ and the random variable $\Y_{\a, c} = \sigma^{(\a)}_{\tau_c}$ with an independent $\tau_c$ as above, we have
$$\frac{c}{c + n^\a} \; =\;\esp [e^{-n \Y_{\a, c}}]$$
for every $n\ge 0.$ We then obtain analogously 
$$G_{\a, c} (x) \; = \; \frac{1}{2\pi c}\;\esp\lcr 1\; +\; \sum_{n\ge 1} (e^{\ii nx} + e^{-\ii nx})\, e^{-n \Y_{\a, c}}\rcr\; = \; \frac{1}{2\pi c}\;\esp\lcr \frac{1- e^{-2\Y_{\a, c}}}{1 - 2 \cos (x)\, e^{-\Y_{\a, c}} + e^{-2\Y_{\a, c}} }\rcr.$$
This alternative representation also shows at once that $G_{\a, c} (x)$ is positive and decreasing on $(0,\pi).$ Moreover, the density of the random variable $\Y_{\a, c} \elaw c^{-1/\a} \Y_{1, \a}$ can be expressed from Formula (4.10.1) in \cite{GKMR} in terms of the Mittag-Leffler function as
$$c\, x^{\a - 1} E_{\a, \a} (-c x^\a)$$
on $\rl^+.$ Therefore,
\begin{eqnarray}
\label{ML}
G_{\a, c} (x) & = & \frac{1}{2\pi}\;\int_0^\infty \frac{1- e^{-2t}}{1 - 2 \cos (x)\, e^{-t} + e^{-2t}} \, t^{\a - 1} E_{\a, \a} (-c t^\a)\, dt\\
\nonumber & = & \frac{1}{2\pi c}\; +\;\frac{1}{\pi} \int_0^\infty \frac{e^{t} \cos(x) - 1}{1 - 2 \cos (x)\, e^{t} + e^{2t}} \, t^{\a - 1} E_{\a, \a} (-c t^\a)\, dt
\end{eqnarray}
and we recover Proposition 3.1 in \cite{LP} - beware that the constant before the integral therein should be $1/\pi$ and not $1/(\pi c)$ because of (3.6) in \cite{LP}, in the case $\a\in (0,1].$ In the case $\a > 1,$ the formula (\ref{ML}) is also true, by Formula (4.10.1) in \cite{GKMR} which implies
$$\frac{1}{c + n^\a}\; =\; \int_0^\infty e^{-n t}\, t^{\a - 1} E_{\a, \a} (-c t^\a)\, dt, \qquad n\ge 1,$$
and a switching between sum and integral which is justified by Fubini's theorem. \\

(b) For $\a\in (0,2),$ the asymptotic expansion (4.9.13) in \cite{GKMR} shows that the integral in (\ref{ML}) is absolutely convergent and can be differentiated  from the inside : we obtain
$$G_{\a, c}' (x) \; = \; \frac{\sin(x)}{\pi}\;\int_0^\infty \frac{e^{-3t}- e^{-t}}{(1 - 2 \cos (x)\, e^{-t} + e^{-2t})^2} \; t^{\a - 1} E_{\a, \a} (-c t^\a)\, dt,$$
which implies $G_{\a,c}'(\pi-) = 0.$ By the aforementioned formula (A.5) in \cite{LP}, we also have $G_{2,c}'(\pi-) = 0.$ This shows that there does not exist any positive Radon measure $\mu$ on $[0,\infty)$ such that 
$$G_{\a,c} (x)\; =\; \int_0^\infty e^{-xt} \, d\mu (t)$$
for $x\in (0,\pi)$ and, by Bernstein's theorem, that $G_{\a,c}(x)$ cannot be extended into a completely monotonic function on the whole half-line.\\

(c) The complete monotonicity of the non-periodic Green function $H_{\a,c}(x)$ on the half-line is certainly folklore, although we could not locate it precisely in the literature. The positivity and monotonicity of $H_{\a,c}(x)$ was obtained in Lemma A.4 (ii) of \cite{FL} as a consequence of the same property for the underlying semigroup. The argument does not work for the further derivatives, because the density of the semigroup is typically flat at zero. The same problem occurs for $G_{\a,c} (x)$, whose complete monotonicity cannot be deduced from (\ref{JTP}) or (\ref{ML}) since the function inside the expectation or the integral is also flat at zero. We stress that the convexity of $H_{\a, c}(x)$ is needed to deduce from (\ref{Sum}) the monotonicity of $G_{\a,c}(x).$ Let us finally refer to \cite{AMV} for the complete monotonicity of another fractional Green function related to the sum of so-called Caputo-Djrbashian fractional derivatives. \\

(d) By multiplicative convolution, the formula (\ref{HC}) shows that for every $\a \in (0,2]$ and $c > 0,$ the function $c H_{\a,c} (x)$ on $\rl$ is the density of the independent product 
$$c^{-1/\a}\,\E\,\times\, \X_\a$$
where the random variable $\E$ is the double exponential with density $(1/2)e^{-\vert x\vert}$ on $\rl,$ the random variable $\X_\a$ has density
$$\frac{2 \sin(\pi\a/2)}{\pi} \lpa \frac{t^{\a-1}}{1 + 2\cos(\pi\a/2)\, t^\a + t^{2\a}}\rpa$$
on $(0,\infty)$ for $\a\in(0,2),$ and $\X_2 \equiv \Un.$ This factorization can also be obtained by Mellin inversion and we leave the details to the reader.\\

(e) For $\a\in (2,4],$ the Laplace transform for the sine function implies
$$\frac{1}{c + u^\a}\; =\; \frac{1}{c} \; \esp\lcr \int_0^\infty e^{-u^2 \sigma^{(\gamma)}_t} \! \sin (ct)\, dt\rcr$$
with $\gamma = \a/4 \in (1/2,1].$ This leads to
$$G_{\a, c} (x) \; =\; \frac{1}{2\pi c}\;\esp\lcr \int_0^\infty \sin (ct) \prod_{n\ge 0} \!\lpa 1 + 2 \cos (x) e^{-(2n+1) \sigma^{(\gamma)}_t}\!\! + e^{-(4n+2) \sigma^{(\gamma)}_t}\rpa\!\lpa 1- e^{-(2n+1) \sigma^{(\gamma)}_t}\rpa dt\rcr.$$  
However, it is not clear whether such a formula can be of any help to solve the open problem stated in \cite{LP} on the number of zeroes of $G_{\a, c} (x)$ for $\a\in (2,4).$
 }

\end{rem}


\begin{thebibliography}{10}  

\bibitem{AAR}
G.~E.~Andrews, R.~Askey and R.~Roy. {\em Special functions.} Cambridge University Press, Cambridge, 1999.

\bibitem{AMV}
V.~V.~Anh and R.~McVinish. Completely monotone property of fractional Green functions. {\em Frac. Calc. Appl. Anal.} {\bf 6} (2), 157-173, 2003. 

\bibitem{FL}
R.~L. Frank and E.~Lenzmann. Uniqueness of non-linear ground states for fractional Laplacian in $\rl.$ {\em Acta Math.} {\bf 210}, 261-318, 2013

\bibitem{GKMR}
R.~Gorenflo, A.~A.~Kilbas, F.~Mainardi and S.~V. Rogosin. {\em Mittag-Leffler functions, related topics and applications.} Springer Verlag, Heidelberg, 2020.

\bibitem{LP}
U.~Le and D.~E.~Pelinovsky. Green's function for the fractional KdV equation on the periodic domain via Mittag-Leffler's function. To appear in {\em Fractional Calculus and Applied Analysis}.

\bibitem{ST}
G.~Samorodnitsky and M.~S.~Taqqu. {\em Stable Non-Gaussian Random Processes.} Chapman \& Hall, New-York, 1994.

\bibitem{Z}
A.~Zygmund. {\em Trigonometrical series}. Chelsea Publishing, New York, 1952.

\end{thebibliography}
\end{document}